# Products of Factorials Modulo $p$

Florian Luca[1]     and     Pantelimon Stănică[2] [*]

[1] Instituto de Matematicas de la UNAM, Campus Morelia
Apartado Postal 61-3 (Xangari) CP 58 089, Morelia, Michoacan, Mexico
[2] Auburn University Montgomery, Department of Mathematics,
Montgomery, AL 36124-4023, USA; e-mail: stanica@sciences.aum.edu

November 3, 2018

## 1  Introduction

Let $p$ be a fixed odd prime and let $s$ and $t$ be fixed positive integers which depend on $p$. Consider the following subset of the elements of $\mathbf{Z}_p^*$

$$P_{s,t}(p) = \{x_1! \cdot x_2! \cdot \cdots \cdot x_t! \pmod{p} \mid x_i \geq 1 \text{ for } i = 1, 2, \ldots, t, \text{ and } \sum_{i=1}^{t} x_i = s\}. \quad (1)$$

The problem that we investigate in this note is loosely the following: given $p$, find sufficient conditions that the parameters $s$ and $t$ should satisfy such as to ensure that $P_{s,t}(p)$ contains the entire $\mathbf{Z}_p^*$.

Let $\varepsilon > 0$ be any small number. Throughout this paper, we denote by $c_1$, $c_2$, ... computable positive constants which are either absolute or depend on $\varepsilon$. From the way we loosely formulated the above problem, it is easily seen that its answer is easily decidable if either both $s$ and $t$ are very small (with respect to $p$) or very large with respect to $p$. For example, if $s < (\log(p))^{2-\varepsilon}$, then it is clear that $P_{s,t}(p)$, or even the union of all $P_{s,t}(p)$ for all allowable values of $t$, cannot possibly contain the entire $\mathbf{Z}_p^*$ when $p$ is large. Indeed, the reason here is that the cardinality of the union of all $P_{s,t}(p)$ for all allowable values of $t$ is at most $p(s) = O(\exp(c_1\sqrt{s}))$ and this is much smaller than $p$ when $p$ is large. Here, we denoted by $p(s)$ the number of unrestricted partitions of $s$. It is also obvious that $P_{s,t}(p)$ does not generate the entire $\mathbf{Z}_p^*$ (for any $s$) when $t = 2$. Moreover, the fact that there exist infinitely many prime numbers $p$ for which the smallest non-quadratic residue modulo $p$ is at least $c_2 \log(p)$, shows that if one wants to generate the entire $\mathbf{Z}_p^*$ out of $P_{s,t}(p)$, then one should allow in (1) partitions of $s$ where $\max(x_i)_{i=1}^t$ is at least $c_2 \log(p)$. In particular, $s$ and $t$ cannot be too close to each other. Indeed, if $p$ is such a prime and the maximum value of the $x_i$'s allowed in (1) is at most $c_2 \log(p)$, then all the numbers in $P_{s,t}(p)$ will be quadratic residues modulo $p$, and in particular $P_{s,t}(p)$ cannot contain the entire $\mathbf{Z}_p^*$. On the other hand, when both $s$ and $t$ are very large, for example $t > 3p$ and $s > p^{3/2+\varepsilon}$, then an immediate argument based on the known upper bounds for the size of the smallest primitive root modulo $p$ shows that $P_{s,t}(p)$ does indeed cover the entire $\mathbf{Z}_p^*$ when $p$ is large. Thus, the question becomes interesting when we search for *small* values of both $s$ and $t$ for which $P_{s,t}(p)$ does cover the entire $\mathbf{Z}_p^*$.

[*]Also Associated to the Institute of Mathematics of Romanian Academy, Bucharest, Romania

This question was inspired by the paper of the second author [8]. In that paper, the problem investigated was the exponent at which a prime number $p$ divides some generalized Catalan numbers. However, the question of whether a certain subset of Catalan numbers, namely the numbers of the form

$$\frac{1}{p} \cdot \binom{p}{m_1, m_2, \ldots, m_t} \tag{2}$$

covers the entire $\mathbf{Z}_p^*$ was not investigated in [8]. Here, the numbers appearing in (2) are all the non-trivial multinomial coefficients. In our notation, this question reduces to whether or not

$$\bigcup_{t\geq 2} P_{p,t}(p) \qquad (3)$$

is the entire $\mathbf{Z}_p^*$. As a byproduct of our results, we show that the set (3) is indeed the entire $\mathbf{Z}_p^*$, for $p \neq 5$.

Our main results are the following:

**Theorem 1.**

*Let $\varepsilon > 0$ be arbitrary. There exists a computable positive constant $p_0(\varepsilon)$ such that whenever $p > p_0(\varepsilon)$, then $P_{s,t}(p) = \mathbf{Z}_p^*$ for all $t$ and $s$ such that $t > p^\varepsilon$ and $s - t > p^{1/2+\varepsilon}$.*

The above result is certainly very far from best possible. We believe that the exponent $1/2$ appearing at the power of $p$ in the lower bound for $s - t$ can be replaced by a much smaller one, or even maybe that the statement of Theorem 1 above remains true when $s - t > p^{2\varepsilon}$. We have not been able to find an argument to prove such a claim.

**Theorem 2.**

*The set (3) is the entire $\mathbf{Z}_p^*$, if $p \neq 5$ is prime.*

The trick in proving Theorem 2 is to detect a small value of $p_0$ such that Theorem 2 holds for $p > p_0$, and then to test the claim for all prime numbers $p$ from 2 up to $p_0$.

Theorem 1 above shows, in particular, that the set (3) (even a very small subset of it) is the entire $\mathbf{Z}_p^*$ when $p$ is large. As an example for Theorem 1, we can easily prove that if 2 is a primitive root modulo $p$, then $A \cup B$, where

$$A = \left\{ 2^u \left(\frac{p-1}{2}\right)! \mid 1 \leq u \leq \frac{p-1}{2} \right\}$$
$$B = \left\{ 2^v \left(\frac{p-3}{2}\right)! \mid 0 \leq v \leq \frac{p-3}{2} \right\}$$

cover the entire $\mathbf{Z}_p^*$. We see first that $A$ and $B$ each contain $\frac{p-1}{2}$ distinct residues modulo $p$. The intersection $A \cap B$ is empty, when 2 is a primitive root modulo $p$. We omit the details. What is interesting is that, in general, we can cover easily all the even residues, and the odd residues from the first half of $\mathbf{Z}_p^*$, since

$$\frac{1}{p}\binom{p}{2, 2k-1, p-2k-1} \equiv k \pmod{p}$$
$$\frac{1}{p}\binom{p}{1, 1, 2k-1, p-2k-1} \equiv 2k \pmod{p},$$

for any $1 \leq k \leq \frac{p-1}{2}$.

Related to our work, we recall that the behaviour of the sequence $n! \pmod p$ was recently investigated in [1].

## 2  The Proofs of the Theorems

The main idea behind the proofs of both Theorems 1 and 2 is to find a suitable list $x_1, x_2, \ldots, x_t$ consisting of many small numbers and each one of them repeated a suitable number of times, such that we can modify (in a sense that will be made precise below) the fixed element given by formula (1) for this list of elements $x_1, x_2, \ldots, x_t$ in enough ways (such that, of course, these modified numbers do not get outside $P_{s,t}(p)$) so that to ensure that in the end, we have obtained all the congruence classes in $\mathbf{Z}_p^*$.

Here is the basic operation by which we can modify a fixed element, call it

$$F = \prod_{i=1}^{t} x_i! \tag{4}$$

in such a way as to obtain, hopefully, new elements in $P_{s,t}(p)$.

(M) *Assume that $i_1 < i_2 < \ldots i_j$ and $l_1 < l_2 < \cdots < l_j$ are two disjoints subsets of indices*

*in* $\{1, 2, \ldots, t\}$. *Then,*

$$\Big(\prod_{s=1}^{j}(x_{l_s}+1)\Big)\Big(\prod_{s=1}^{j}x_{i_s}\Big)^{-1} \cdot F = x_1! \cdots (x_{l_1}+1)! \cdots (x_{i_1}-1)! \cdots \\ (x_{l_2}+1)! \cdots (x_{i_2}-1)! \cdots x_t! = F' \in P_{s,t}(p). \quad (5)$$

In general, we shall always apply formula (5) with $x_{l_1} = \cdots = x_{l_j} = 1$. With this convention, we may eliminate the initial number $F$, take inverses in (5) above, and then reformulate the question as follows:

**Question:** *Is it true that for suitable integers $t$ and $s$ (satisfying, for example, the hypothesis of Theorem 1) we can find some positive integers $x_1, x_2, \ldots, x_t$ summing to $s$, such that every non-zero residue class modulo $p$ can be represented by a number of the form*

$$\prod_{r=1}^{j}\Big(\frac{x_{i_r}}{2}\Big) \quad (6)$$

*where a subset of indices $\{i_1, i_2, \ldots, i_j\}$ of $\{1, 2, \ldots, t\}$ in (6) can be any subset as long as there exists another subset of $j$ indices $\{l_1, l_2, \ldots, l_j\}$ disjoint from $\{i_1, i_2, \ldots, i_j\}$ for which $x_{l_r} = 1$ for all $r = 1, 2, \ldots, j$?*

**The Proof of Theorem 1.** All we have to show is that if the parameters $s$ and $t$ satisfy the hypothesis of Theorem 1, then we can construct a list of elements $x_1, x_2, \ldots, x_t$ for which the answer to the above question is affirmative. Fix $\varepsilon > 0$ and a positive integer $k$ with $\frac{1}{k} < \varepsilon < \frac{2}{k}$. From now on, all positive constants $c_1, c_2, \ldots,$ which will appear will be computable and will depend only on $k$. We shall show that if $p$ is large enough with respect to $k$, then we can construct a good sublist of numbers $x_1, x_2, \ldots, x_t$ in the following manner:

1. *We first take and repeat exactly two times each of the prime numbers $x_i$ up to $p^{1/k}$.*

2. *We then adjoin some even numbers $x_j$, each one of them smaller than $p^{1/2+1/k}$ but such that the totality of those (counted with multiplicities) does not exceed $c_1 \log(\log(p))$.*

3. *The numbers of the form (6), where the $x_i$'s are from the lists 1 and 2 and the maximum length $j$ of a product in (6) is not more than $2k + 2c_1 \log(\log(p))$ cover the entire $\mathbf{Z}_p^*$.*

It is clear that if we can prove the existence of a list satisfying 1-3 above, then we are done. Indeed, we may first adjoin at the sublist consisting of the numbers appearing at 1 and 2 above a number of about $2k + 2c_1 \log(\log(p))$ values of $x_i$ all of them equal to 1. The totality of all these numbers (the ones from 1, 2 and these new values of $x_i$ all equal to 1) counted with their multiplicities, so far, is certainly not more than

$$c_2 \frac{p^{1/k}}{\log(p)} + 2k + 4c_1 \log(\log(p)) < p^{\varepsilon} - 1 < t - 1, \tag{7}$$

while their sum is at most

$$c_3 \frac{p^{2/k}}{\log(p)} + 2k + 2c_1 \log(\log(p)) + 2c_1 \log(\log(p)) p^{1/2+1/k} < p^{1/2+\varepsilon} - 1 < s - t - 1, \tag{8}$$

for large $p$. At this step, we may finally complete the above list with several other values of the $x_i$ equal to 1 until we get a list with precisely $t-1$ numbers, which is possible by inequality (7) above, and set the last number of the list to be equal to

$$x_t = s - \sum_{i=1}^{t-1} x_i,$$

which is still positive by inequality (8) above.

To show the existence of a sublist with properties 1-3 above, we start with the set

$$A = \{n \mid n < p^{1/k} \text{ and } n \text{ prime}\}. \tag{9}$$

The numbers from $A$ will form the sublist mentioned at 1 above but, so far, we take each one of them exactly once. Let

$$B_1 = \left\{ \frac{n_1}{2} \cdot \frac{n_2}{2} \cdot \ldots \cdot \frac{n_k}{2} \;\middle|\; n_i \in A,\; n_i \neq n_j \text{ for } 1 \leq i \neq j < k \right\}. \tag{10}$$

We first notice that each value of $n \in A$ appears at most $k$ times in an arbitrary product in $B_1$. We now show that $b_1 = \#B_1$ is large. Indeed, the set $B_1$ will certainly contain all the numbers of the form

$$\frac{p_1}{2} \cdot \frac{p_2}{2} \cdot \ldots \cdot \frac{p_k}{2} = 2^{-k} \cdot p_1 \cdot p_2 \cdot \ldots \cdot p_k, \tag{11}$$

where $p_i$ is an arbitrary prime subject to the condition

$$p_i \in \left( \frac{p^{1/k}}{2^i}, \frac{p^{1/k}}{2^{i-1}} \right) \quad \text{for } i = 1,\, 2,\, \ldots,\, k. \tag{12}$$

Moreover, notice that the residue classes modulo $p$ of the elements of the form (11), where the primes $p_i$ satisfy conditions (12), are all distinct. Indeed, the point is that if two of the numbers of the form (11) coincide modulo $p$, then, after cancelling the $2^{-k}$, we get two residue classes of integers which coincide modulo $p$. Now each one of these two integers is smaller than $p$, therefore if they coincide modulo $p$, then they must be, in fact, equal. Now the fact that they are all distinct follows from the fact that their prime divisors $p_i$ satisfy condition (12). Applying the Prime Number Theorem to estimate from below the number of primes in each one of the intervals appearing in formula (12), we get

$$b_1 > c_4 \frac{p}{(\log(p))^k} > \frac{p}{(\log(p))^{k+1}}, \tag{13}$$

whenever $p > c_5$. If $B_1$ is the entire $\mathbf{Z}_p^*$, then we are done. Assume that it is not so. We construct recursively a (finite) increasing sequence of subsets $B_m$ for $m \geq 1$ in the following way:

Assume that $B_m$ has been constructed and set $b_m = \#B_m$. Assume that $b_m < p-1$ (that is, $B_m$ is not the entire $\mathbf{Z}_p^*$ already). Then, we have the following trichotomy:

i. If $b_m \geq p/2$, then set $B_{m+1} = B_m \cdot B_m$, and notice that $B_{m+1} = \mathbf{Z}_p^*$ and we can no longer continue.

ii. If $b_m < p/2$ and there exists an even number $a < p^{1/2+1/k}$ such that $a/2 \notin B_m \cdot B_m^{-1}$, then set $a_m = a$, add $a$ to the list of the $x_i$'s (as one of the numbers from 2. above) and let

$$B_{m+1} = B_m \cup \frac{a_m}{2} \cdot B_m. \tag{14}$$

Notice that

$$b_{m+1} \geq \min(p-1,\ 2b_m). \tag{15}$$

iii. If $b_m < p/2$ and all even numbers $a$ up to $p^{1/2+1/k}$ have the property that $a/2$ is already in $B_m \cdot B_m^{-1}$, we choose the even number $a$ smaller than $p^{1/2+1/k}$ for which the number of representations of $a/2$ of the form $x \cdot y^{-1}$ with $x,\ y \in B_m$ is minimal. Then, we set $a_m = a$, add $a$ to the list of the $x_i$'s (as one of the numbers from 2. above), set

$$B_{m+1} = B_m \cup \frac{a_m}{2} \cdot B_m, \tag{16}$$

and notice that

$$b_{m+1} \geq \min\left(p-1,\ \frac{4b_m}{3}\right). \tag{17}$$

In i-iii above we have used the set-theoretic notation, namely that if $U$ and $V$ are two subsets of $\mathbf{Z}_p^*$, we have denoted by $U \cdot V$ the set of all elements of $\mathbf{Z}_p^*$ of the form $u \cdot v$ with $u \in U$ and $v \in V$, and by $U^{-1}$ the set of all elements of the form $u^{-1}$ for $u \in U$.

We have to justify that i-iii above do indeed hold. Notice that i and ii are obvious. The only detail we have to justify is that inequality (17) indeed holds in situation iii. For this, we use the following Lemma due to Sárközy (see [7]):

**Lemma S.**

Let $p$ be a prime number, $u$, $v$, $S$, $T$ be integers with $1 \leq u$, $v \leq p-1$, $1 \leq T \leq p$, furthermore $c_1$, $c_2$, ..., $c_u$ and $d_1$, $d_2$, ..., $d_v$ are integers with

$$c_i \not\equiv c_j \pmod{p}, \quad \text{for } 1 \leq i < j \leq u, \tag{18}$$

and

$$d_i \not\equiv d_j \pmod{p}, \quad \text{for } 1 \leq i < j \leq v. \tag{19}$$

For any integer $n$, let $f(n)$ denote the number of solutions of

$$c_x \cdot d_y \equiv n \pmod{p}, \quad 1 \leq x \leq u, \ 1 \leq y \leq v. \tag{20}$$

Then,

$$\left| \sum_{n=S+1}^{S+T} f(n) - \frac{uvT}{p} \right| < 2(puv)^{1/2} \log(p). \tag{21}$$

We apply Lemma S above with $u = v = b_m$, $c_1$, $c_2$, ..., $c_u$ all the residue classes in $B_m$ and $d_1$, $d_2$, ..., $d_u$ all the residue classes in $B_m^{-1}$. We also set $S = 0$ and $T$ to be the largest integer smaller than $p^{1/2+1/k}/2$. Clearly, $T > p^{1/2+1/k}/3$. Since we are discussing situation iii above, we certainly have $f(n) \geq 1$ for all positive integers $n$ up to $T$. Let $M = \min(f(n) \mid 1 \leq n \leq T)$ and then $a_m = 2c$ where $f(c) = M$. Denote $b_m$ by $b$. We apply inequality (21) to get

$$M < \frac{b^2}{p} + \frac{2b\sqrt{p}\log(p)}{T}. \tag{22}$$

We first show that

$$\frac{2b\sqrt{p}\log(p)}{T} < \frac{b^2}{3p} \qquad (23)$$

holds. Indeed, since $T > p^{1/2+1/k}/3$ and $b = b_m \geq b_1 > \frac{p}{(\log(p))^{k+1}}$ (by inequality (13)), it follows that in order for (23) to hold it suffices that

$$54(\log(p))^{k+2} < p^{1/k}, \qquad (24)$$

which is certainly satisfied when $p > c_6$. Thus, inequalities (22) and (23) show that

$$M < \frac{4b^2}{3p} < \frac{2b}{3}, \qquad (25)$$

where the last inequality in (25) follows because $b < p/2$. In particular,

$$b_{m+1} = \#(B_m \cup c \cdot B_m) \geq b_m + (b_m - M) \geq 2b - \frac{2b}{3} = \frac{4b}{3}, \qquad (26)$$

which proves inequality (17).

The combination of (13), (15) and (17) show that

$$b_{m+1} > \left(\frac{4}{3}\right)^m b_1 > \left(\frac{4}{3}\right)^m \frac{p}{(\log(p))^{k+1}} \qquad (27)$$

holds as long as $b_m < p/2$. Now notice that the inequality

$$\left(\frac{4}{3}\right)^m > \frac{(\log(p))^{k+1}}{2} \qquad (28)$$

will happen provided that $m > c_7 \log(\log(p))$, where one can take $c_7 = \frac{k+1}{\log(4/3)}$, for example, and for such large $m$ inequality (27) shows that $b_{m+1} > p/2$. In particular, situations ii or iii above will not occur more than $c_7 \log(\log(p))$ times after which we arrive at a point where we apply situation i to construct $B_{m+1}$ and we are done. Clearly, i-iii and the above arguments prove the existence of a sublist of the $x_i$'s satisfying conditions 1-3, which finishes the proof of Theorem 1.

**The Proof of Theorem 2.** We follow the method outlined in the proof of Theorem 1. Thus, it suffices to find a list of positive integers, say $A := \{x_1, x_2, \ldots, x_s\}$, with

$$U := \sum_{i=1}^{s} x_i < p, \tag{29}$$

and such that for every $m \in \mathbf{Z}_p^*$ there exists a subset $I \subseteq \{1, 2, \ldots, s\}$ for which

$$m \equiv \prod_{i \in I} x_i! \pmod{p}. \tag{30}$$

**Step 1.** We start with a set $A_1$ of distinct positive integers such that

$$U_1 := \sum_{x \in A_1} x \tag{31}$$

is not too large, and set $B_1 := A_1 \cdot A_1 \pmod{p}$. For $m \geq 1$, we construct inductively the sets $A_m$ and $B_m$ by the method explained in the proof of Theorem 1. We set $b_m := \#B_m$, $s_m := b_m/p$, and we choose the parameter $T$ to be of the form

$$T := 2\lfloor \lambda \sqrt{p} \log p \rfloor + 1, \tag{32}$$

where $\lambda > 2$ is some parameter, for which we shall specify later an optimal value, and $\lfloor x \rfloor$ is the *floor function of $x$*, that is the largest integer which is less than or equal to $x$. From the way the sets $A_m$ and $B_m$ are constructed for $m \geq 1$, it follows that as long as $s_m < 1/2$, $A_{m+1}$ is obtained from $A_m$ by adjoining to it just one element $a_m$ of size no larger than $T$, and then $B_{m+1}$ is taken to be $B_m \cup a_m \cdot B_m \pmod{p}$. Thus,

$$U_{m+1} := \sum_{x \in A_{m+1}} x \leq T + \sum_{x \in A_m} x = T + U_m, \quad \text{for } m \geq 1, \tag{33}$$

and therefore

$$U_{m+1} \leq mT + U_1, \qquad (34)$$

and the above inequality (34) holds for all $m \geq 1$ as long as $s_m < 1/2$. However, by formula (22) and our choice for $T$, it follows that when constructing $A_{m+1}$ out of $A_m$, we choose the parameter $M$ in such a way that

$$M < \frac{b_m^2}{p} + \frac{2b_m\sqrt{p}\log p}{T} < b_m\left(s_m + \frac{1}{\lambda}\right),$$

therefore inequality (26) now shows that

$$b_{m+1} \geq 2b_m - M > b_m\left(\left(2 - \frac{1}{\lambda}\right) - s_m\right).$$

Hence,

$$s_{m+1} > (\beta - s_m)s_m, \qquad (35)$$

where

$$\beta := \beta(\lambda) := 2 - \frac{1}{\lambda} = \frac{2\lambda - 1}{\lambda}. \qquad (36)$$

Of course, the above construction will be repeated only as long as $s_m < 1/2$. If we denote by $n$ the largest positive integer such that $s_n < 1/2$, then $s_{n+1} \geq 1/2$, therefore the last set $B_{n+2}$, which is the entire $\mathbf{Z}_p^*$, is taken to be $B_{n+1} \cdot B_{n+1} \pmod{p}$, i.e., $A_{n+2}$ is taken to be the list of elements $A_{n+1}$, but now each one of them is repeated twice. Thus,

$$U_{n+2} \leq 2U_{n+1} \leq 2(nT + U_1). \qquad (37)$$

From these arguments, it follows that in order to insure that $U_{n+2}$ is not larger than $p-1$, it suffices to check that

$$2(nT + U_1) < p. \qquad (38)$$

The number $U_1$ can be easily computed in terms of $A_1$, therefore all we need in order to check that inequality (38) holds, is a good upper bound on $n$ in terms of $A_1$. We recall that $n$ is the largest positive integer with $s_n < 1/2$, where the sequence $(s_m)_{m \geq 1}$ has initial term $s_1 := b_1/p$ and satisfies the recurrence (35).

**Step 2.** We give an upper bound on $n$. Since $\lambda > 2$, it follows that $\beta > 3/2$, therefore inequality (35) shows that $s_{m+1} > s_m$ as long as $s_m < 1/2$. By (35), we also have

$$s_{k+1} > \beta^k\left(1 - \frac{s_k}{\beta}\right), \qquad \text{for } k = 1, 2, \ldots, n,$$

therefore

$$s_{n+1} > \beta^n s_1 \prod_{k=1}^{n}\left(1 - \frac{s_k}{\beta}\right). \tag{39}$$

Since $s_k < 1/2$ for $k = 1, 2, \ldots, n$, it follows that

$$\frac{s_k}{\beta} < \frac{1}{2\beta} = \frac{\lambda}{2(2\lambda - 1)}. \tag{40}$$

The inequality

$$1 - x > e^{-\mu x} \tag{41}$$

holds for all $x$ in the interval $\left(0, \frac{\lambda}{2(2\lambda - 1)}\right)$ with some value $\mu := \mu(\lambda)$, and the best value of $\mu$ is precisely

$$\mu := -\frac{\log(1-x)}{x}\bigg|_{x=\frac{1}{2\beta}} = \frac{2(2\lambda - 1)}{\lambda} \cdot \log\left(\frac{4\lambda - 2}{3\lambda - 2}\right). \tag{42}$$

The fact that the best value of $\mu$ for which inequality (41) holds with all $x$ in the interval $\left(0, \frac{1}{2\beta}\right)$ is indeed the one given by formula (42) follows from the fact that the function $x \to -\frac{\log(1-x)}{x}$ is decreasing in the interval $\left(0, \frac{1}{2\beta}\right]$. Thus,

$$\log s_{n+1} > n\log\beta + \log s_1 + \sum_{k=1}^{n}\log\left(1 - \frac{s_k}{\beta}\right) > n\log\beta + \log s_1 - \frac{\mu}{\beta}\sum_{k=1}^{n} s_k. \tag{43}$$

We now find an upper bound on the sum appearing in the right hand side of inequality (43). Notice that since $\lambda > 1/2$, it follows that whenever $s_m < 1/2$, one also has

$$s_{m+1} > (\beta - s_m)s_m > (1+\rho)s_m, \tag{44}$$

where the best $\rho := \rho(\lambda)$ is given by

$$\beta - \frac{1}{2} = 1 + \rho,$$

or, equivalently,

$$\rho := \beta - \frac{3}{2} = \frac{1}{2} - \frac{1}{\lambda} = \frac{\lambda - 2}{2\lambda},$$

and

$$1 + \rho = \frac{3\lambda - 2}{2\lambda}. \tag{45}$$

In particular,

$$s_{n-1} < \frac{1}{1+\rho} s_n$$

holds, and if $k$ is any positive integer less than $n$, then

$$s_{n-k} < \left(\frac{1}{1+\rho}\right)^k s_n$$

holds. Thus,

$$\sum_{k=1}^{n} s_k < s_n \sum_{k \geq 0} \left(\frac{1}{1+\rho}\right)^k < \frac{1}{2}\frac{\rho+1}{\rho} = \frac{3\lambda - 2}{2(\lambda - 2)}.$$

The above calculations show that

$$\log a_{n+1} > n \log \beta + \log s_1 - \mu \cdot \frac{(3\lambda - 2)\lambda}{2(2\lambda - 1)(\lambda - 2)} = n \log \beta + \log s_1 - \gamma, \tag{46}$$

where

$$\gamma := \gamma(\lambda) := \mu \cdot \frac{(3\lambda - 2)\lambda}{2(2\lambda - 1)(\lambda - 2)} = \frac{(3\lambda - 2)}{(\lambda - 2)} \cdot \log\left(\frac{4\lambda - 2}{3\lambda - 2}\right). \tag{47}$$

Thus, if we choose $n$ such that

$$n \log \beta + \log s_1 - \gamma \geq \log(1/2), \tag{48}$$

then we are sure that $s_{n+1} > 1$. Inequality (48) is equivalent to

$$n \log \beta > -\log(2s_1) + \gamma,$$

hence, to

$$n > \frac{1}{\log \beta}\left(-\log(2s_1) + \gamma\right). \tag{49}$$

Thus, we may write

$$n_0 := 1 + \left\lfloor \frac{1}{\log \beta}\left(-\log(2s_1) + \gamma\right) \right\rfloor, \tag{50}$$

and conclude that $n \leq n_0$. Thus, inequality (38) will be satisfied provided that

$$n_0 T + U_1 < \frac{p}{2} \tag{51}$$

holds, where $n_0$ is given by formula (50).

**Step 3.** Here, we show that we can do the above construction for $p > 3.242 \cdot 10^6$. All we need to do is to explain how we choose $A_1$, to give a lower bound on $s_1$ and an upper bound on $U_1$, and to check that inequality (51) holds. From here on, we write $x := p$ and $y := \sqrt{\frac{x}{2}}$, and we assume that $x > 3.242 \cdot 10^6$, therefore that $y > 1163$. We choose

$$A_1 := \{q \mid q \text{ is prime and } q \leq y\}, \tag{52}$$

and therefore

$$B_1 := \{q_1 q_2 \mid q_1 < q_2 \text{ and } q_1, q_2 \in A\}.$$

It is clear that the elements of $B_1$ are in distinct congruence classes in $\mathbf{Z}_p^*$, therefore we may consider $B_1$ as a subset of $\mathbf{Z}_p^*$ and its cardinality is precisely

$$b_1 := \binom{\pi(y)}{2} = \frac{\pi(y)(\pi(y) - 1)}{2} > \frac{(\pi(y) - 1)^2}{2}, \tag{53}$$

where $\pi(y)$ is the number of primes up to $y$. We first show that

$$\pi(y) - 1 > \frac{y}{\log y - 0.5}. \tag{54}$$

We recall that Rosser and Schoenfeld [6] showed a long time ago that both inequalities

$$\pi(x) > \frac{x}{\log x - 0.5}, \quad \text{for all } x \geq 67 \tag{55}$$

and

$$\pi(x) < \frac{x}{\log x - 1.5}, \quad \text{for all } x > e^{3/2} \tag{56}$$

hold, but the above inequalities were slightly strengthened recently by Panaitopol (see [5]) who showed that in fact both inequalities

$$\pi(x) > \frac{x}{\log x - 1 + (\log x)^{-0.5}}, \quad \text{for all } x \geq 59 \tag{57}$$

and

$$\pi(x) < \frac{x}{\log x - 1 - (\log x)^{-0.5}}, \quad \text{for all } x \geq 6 \tag{58}$$

hold. Since $y > 1163 > 59$, we may use inequality (57) and infer that in order for inequality (54) to hold it suffices to check that

$$\frac{y}{\log y - 1 + (\log y)^{-0.5}} > \frac{y}{\log y - 0.5} + 1$$

holds. After some manipulations, the last inequality above is seen to be equivalent to

$$y(0.5 - (\log y)^{-0.5}) > (\log y - 1 + (\log y)^{-0.5})(\log y - 0.5). \tag{59}$$

Since $\exp(2.5^2) < 519 < y$, it follows that $\log y > 2.5^2$, therefore

$$0.5 - (\log y)^{-0.5} > \frac{1}{2} - \frac{1}{2.5} = \frac{1}{10},$$

thus, in order for inequality (59) to hold it suffices that

$$\frac{y}{10} > (\log y - 1 + (2.5)^{-1})(\log y - 0.5),$$

and this last inequality is satisfied whenever $y > 246$. By inequalities (53) and (54), it follows that

$$b_1 > \frac{y^2}{2(\log y - 0.5)^2} = \frac{x}{4(\log \sqrt{x/2} - 0.5)^2} = \frac{x}{\log(x/c_1)^2}, \tag{60}$$

where $c_1 := 2e$. We also notice that in the above computations we only needed that $y > 519$, or that $x = 2y^2 > 538722$. Thus, (60) shows that

$$s_1 = b_1/x > \frac{1}{\log(x/c_1)^2}. \tag{61}$$

We next give an upper bound on $U_1$. We claim that

$$U_1 = \sum_{q \in A_1} q < \frac{y^2}{21.4}. \tag{62}$$

Let $N := \pi(y)$. From an inequality in [4], we know that

$$p_m < m\left(\log m + \log \log m - 1 + \frac{1.8 \log \log m}{\log m}\right), \quad \text{for all } m \geq 13. \tag{63}$$

Here $p_m$ denotes the $m$th prime number. Since $y > 1163$ and $\pi(y) \geq \pi(1163) = 192$, it follows that

$$\log N + \log \log N - 1 + \frac{1.8 \log \log N}{\log N} \geq \log 192 + \log \log 192 - 1 + \frac{1.8 \log \log 192}{\log 192} \sim 6.485 > 6$$

and

$$p_m < 6m \quad \text{for } m = 1, 2, \ldots, 13.$$

Thus,

$$p_m < m(\log N + \log \log N - 1 + \frac{1.8 \log \log N}{\log N})$$

holds for all $m = 1, 2, \ldots, N$. Thus,

$$\sum_{q \in A_1} q = \sum_{m=1}^{N} p_m < (\log N + \log \log N - 1 + \frac{1.8 \log \log N}{\log N}) \sum_{m=1}^{N} m$$
$$= \frac{N(N-1)}{2} \cdot \left( \log N + \log \log N - 1 + \frac{1.8 \log \log N}{\log N} \right) \quad (64)$$
$$= \frac{1}{2} \cdot \pi(y)(\pi(y) - 1) \left( \log \pi(y) + \log \log \pi(y) - 1 + \frac{1.8 \log \log \pi(y)}{\log \pi(y)} \right).$$

Thus, in order to check that (62) holds it suffices to check that

$$\pi(y)(\pi(y) - 1) \left( \log \pi(y) + \log \log \pi(y) - 1 + \frac{1.8 \log \log \pi(y)}{\log \pi(y)} \right) < \frac{y^2}{10.7}. \quad (65)$$

Instead of using Panaitopol's inequality (58), we will use another result which belongs to Dusart [3], stating that for $x > 598$ (the upper bound holds for $x > 1$) we have

$$\frac{x}{\log x} \left( 1 + \frac{0.992}{\log x} \right) \leq \pi(x) \leq \frac{x}{\log x} \left( 1 + \frac{1.2762}{\log x} \right). \quad (66)$$

So it suffices to check that (65) holds when $\pi(y)$ is replaced by $\frac{y}{\log y} \left( 1 + \frac{1.2762}{\log y} \right)$, and we checked that this last inequality is true with the starting value for $y := 970 < 1163$; hence, for $x > 2 \cdot (970)^2 = 1881800$. Notice that inequality (62) simply asserts that

$$U_1 < \frac{x}{42.8}. \quad (67)$$

With the lower bound on $s_1$ given by (61) and the upper bound on $U_1$ given by (67), it follows that inequality (51) will hold provided that

$$\left(1 + \frac{1}{\log \beta}\left(-\log 2 + 2\log\log(x/c_1) + \gamma\right)\right)(2\lambda\sqrt{x}\log x + 1) < \frac{x}{2} - \frac{x}{42.8} = \frac{10.2x}{21.4}, \quad (68)$$

with some $\lambda > 2$, where $\beta$ and $\gamma$ are given in terms of $\lambda$ by formulae (36) and (47), respectively. We did some experiments with Mathematica[1], and the best lower bound on $x$ for which inequality (68) holds was found at $\lambda := 3$, for which inequality (68) is satisfied whenever $x \geq 9.1 \cdot 10^6$. Thus, from now on we assume that $x < 9.1 \cdot 10^6$. To cut the range down from $9.1 \cdot 10^6$ to $3.242 \cdot 10^6$, we proceeded as follows. Assume that $x > 3.242 \cdot 10^6$. Notice that by inequality (51) and the upper bound (67) on $U_1$, it suffices to have

$$nT < \frac{x}{2} - \frac{x}{42.8} = \frac{10.2x}{21.4} \quad (69)$$

holds, where $n$ is the largest index for which $s_n < 1/2$. Since we now have a starting value on $x$, namely $x < 9.1 \cdot 10^6$, it follows, by inequality (61), that

$$s_1 > \frac{1}{\log(9.1 \cdot 10^6/c_1)^2} > \frac{1}{206}. \quad (70)$$

For each $\lambda > 2$, let $n(\lambda)$ be the largest value of $n$ for which $s_n < 1/2$, where the sequence $(s_m)_{m \geq 1}$ now has $s_1 := 1/206$, and satisfies the recurrence relation

$$s_{m+1} = (\beta - s_m)s_m, \quad \text{for all } m \geq 1. \quad (71)$$

Since

$$1 + 2\lambda\sqrt{x}\log x < 2(\lambda + 0.001)\sqrt{x}\log x,$$

---
[1] A Trademark of Wolfram Research

(this is simply because $\lambda > 2$ and $x$ is large), it follows that

$$nT < 2n(\lambda)(\lambda + 0.001)\sqrt{x}\log x,$$

and therefore inequality (69) will be satisfied provided that

$$2n(\lambda)(\lambda + 0.001)\sqrt{x}\log x < \frac{10.2x}{21.4},$$

or, equivalently,

$$\left(\frac{42.8}{10.2}\delta \log x\right)^2 < x, \tag{72}$$

where

$$\delta := \delta(\lambda) := n(\lambda)(\lambda + 0.001). \tag{73}$$

We see that in order for (72) to hold starting with a relatively small value of $x$, we need that the expression $\delta$ be small. We have let $\lambda$ take all the values of the form $2 + i/10$ for $i = 1, 2, \ldots, 30$, and for each one of these values of $\lambda$ we have computed $n(\lambda)$. Out of all these values obtained in this way, we selected the $\lambda$ for which $\delta$ is minimal. The minimal value of $\delta$ found is smaller than 28.62, and by replacing $\delta$ by 28.62 in (72), we got an inequality which holds for all $x \geq 3.242 \cdot 10^6$. Thus, it only remains to check the values of $p$ which are less than $3.242 \cdot 10^6$.

## 3 The computer verification

It suffices to check that for all prime numbers $5 < p < 3.242 \cdot 10^6$, the set

$$\{\prod_{i=1}^{t} m_i! \mid \sum_{i=1}^{t} m_i = p - 1\} \tag{74}$$

covers the entire $\mathbf{Z}_p^*$. Here is a trick that worked for the primes $p$ which are large enough (for example, $p > 6 \cdot 10^3$).

**Lemma.**

*Assume that $a > 1$ is a primitive root modulo $p$ and assume that $v$ and $b$ are positive integers in the interval $(1,\ p-1)$ such that $b \equiv a^v \pmod{p}$ and*

$$v^2 a < p(v - b). \tag{75}$$

*Then, the set given by (74) covers $\mathbf{Z}_p^*$.*

**Proof.** Take $w := \lfloor (p-1)/v \rfloor$, $t := (v-1) + w$, and $m_i := a$ for $i = 1,\ 2,\ \ldots,\ v-1$, and $m_i = b$ for $i = v,\ v+1,\ \ldots,\ t$. Notice first that

$$\sum_{i=1}^{t} m_i = (v-1)a + wb < va + \frac{p}{v}b < p,$$

where the last inequality from the right above follows from (75). Thus, we may complete the $t$-tuple $(m_1,\ \ldots,\ m_t)$ with ones until we get a longer vector whose sum of the coordinates is equal to $p - 1$. Notice also that for each pair of non-negative integers $(\lambda,\ \mu)$ with $\lambda \leq v - 1$ and $\mu \leq w$ we have

$$(a!)^{v-1}(b!)^w = a^\lambda b^\mu \big((a-1)!^\lambda (b-1)!^\mu a!^r b!^s\big),$$

where $r = v - 1 - \lambda$ and $s = w - \mu$. Thus, it suffices to show that every congruence class in $\mathbf{Z}_p^*$ can be represented under the form $a^\lambda b^\mu$ for some non-negative $\lambda$ and $\mu$ with $\lambda \leq v - 1$ and $\mu \leq w$. But clearly, every congruence class in $\mathbf{Z}_p^*$ is of the form $a^t$ for some $t$ in the interval $[1,\ p-1]$ (because $a$ is a primitive root modulo $p$). We may now apply the division with remainder theorem to write

$$t = \mu v + \lambda,$$

where $\lambda \leq v - 1$, and $\mu$ is the integer part of $t/v$. Thus, $\mu \leq w$, and

$$a^t = a^{\mu v + \lambda} = a^\lambda (a^v)^\mu = a^\lambda b^\mu,$$

and the Lemma is therefore proved.

It is not even clear to us that for a given prime $p$ there should exist a primitive root $a$ modulo $p$ and a value of the positive integer $v$ such that inequality (75) is satisfied, although under the GRH, we know that there exist small primitive roots modulo $p$, and recent results on the distribution of $a^v \pmod{p}$ for small $v$ (see [2]) might imply that one can find choices for $a$ and $v$ satisfying (75), if $p$ is large enough.

However, we are not interested in whether or not we can prove that one can find choices for $a$ and $v$ satisfying (75), we rather want to check computationally that this is indeed so for $6 \cdot 10^3 < p < 3.242 \cdot 10^6$. For this, fix a prime $p$. We took the first 25 odd primes and we checked each one of them against being a primitive root modulo $p$. It is clear that at least one of these numbers will be a primitive root modulo $p$ for most of the primes $p$ in our range. We collected all these primes (which are primitive roots modulo $p$) in a set which we called $A$. Now we tried to find a value for $v$. We could have looped over all possible values of $v$, but this would have resulted in a cycle of length $p - 1$ for each $p$, and the computation would have been infeasible. Instead, let $v_0$ be an initial value of $v$ and set $b \equiv a^{v_0} \pmod{p}$. If $v_0$ is good, we are done. If not, we set the next $v$ to be such that

$$v := v_0 + 1 + \left\lfloor \frac{\log p/b}{\log a} \right\rfloor. \tag{76}$$

In a sense, the $v$ shown above is the smallest $v > v_0$ one can choose for which there is a chance for $a^v = a^{v_0} a^{v-v_0} = b a^{v-v_0}$ to be small modulo $p$. We kept on doing this for about $3\sqrt{p}$ times for each $a \in A$. If no good values of $a$ and $v$ were found by this code, then we had the program tell us that $p$ is a "bad" prime. The computation was done with $v_0 = \lfloor \log p / \log a \rfloor$ but a different choice of $v_0$ might give better results.

Now, $\pi(3.242 \cdot 10^6) = 233053$. We list here the 112 "bad" primes between $100^{th}$ and $233053^{th}$ prime, obtained after the first run of the algorithm:

541, 601, 661, 709, 853, 911, 1009, 1021, 1091, 1117, 1171, 1297, 1303, 1399, 1429, 1453, 1531, 1549, 1621, 1811, 2029, 2351, 2383, 2441, 3001, 3299, 3319, 3559, 3709, 3877, 4129, 5749, 5881, 7591, 23911, 31771, 46861, 71761, 71761, 93481, 93481, 103091, 190321, 266701, 267901, 290041, 412387, 448141, 453181, 494101, 509389, 513991, 609757, 661093, 674701, 690541, 698491, 775861, 776179, 781051, 790861, 975493, 1026061, 1035829, 1067557, 1152421, 1162951, 1242361, 1308421, 1309699, 1364731, 1418551, 1444873, 1445137, 1506121, 1520851, 1732669, 1853461, 1863541, 1895011, 1897561, 2057701, 2080597, 2100121, 2149351, 2165671, 2171311, 2175109, 2183833, 2238661, 2248171, 2270641, 2273431, 2312311, 2319241, 2370889, 2441041, 2447761, 2480479, 2535331, 2561731, 2656351, 2697301, 2708581, 2728261, 2800141, 2831011, 2857951, 2868139, 3014371, 3026971, 3126061.

We increased the range for $v$ to $10\sqrt{p}$ and we sieved the above list. The output was now only the bad primes up to (and including) 7591. Then we increased the range of primes, which are probable primitive roots modulo $p$. The list shortened to 25 "bad" primes, namely 541, 601, 661, 709, 853, 911, 1009, 1021, 1091, 1117, 1171, 1297, 1399, 1429, 1453, 1531, 1549, 1621, 1811, 2029, 2351, 2441, 3001, 3319, 5749. These primes were handled by a different method: we wrote a Mathematica program which showed that the union of the sets (for the possible values of $s$)

$$A_p(s) = \left\{ 2^u \left( \frac{p-2s-1}{2} \right)! \,|\, 0 \leq u \leq \left\lfloor \frac{p+2s+1}{4} \right\rfloor \right\} \tag{77}$$

covers the entire $\mathbf{Z}_p^*$, for any $p$ in the remaining set of "bad" primes. In fact, the above sets were shown to cover $Z_p^*$ for all primes up to 1000, except for $p = 5$. We conjecture that the union of (77), for all the possible values of $s$, covers $Z_p^*$ for any prime $p \neq 5$.

**Ackowledgements.** We thank Professor William O. Nowell for help with the programming in C++ for double-checking our computation by Mathematica.